\documentclass[12pt,a4j]{article}

\usepackage[dvips]{graphicx, color}
\usepackage{amsmath, amssymb}
\usepackage{array}
\usepackage{float}

\newtheorem{thm}{Theorem}[section]
\newtheorem{prop}[thm]{Proposition}
\newtheorem{cor}[thm]{Corollary}

\newtheorem{rei}[thm]{Example}
\newtheorem{chu}[thm]{Remark}

\title{{\bf Existence of a regular unimodular triangulation of the edge polytopes of finite graphs}} 
\author{GINJI HAMANO}
\date{}

\begin{document}
\maketitle

\begin{abstract}
In this paper, we give several criteria for the edge polytope of a graph to possess a regular unimodular triangulation in terms of some simple data of the graph. 
We further apply our criteria to several examples of graphs and show that their edge polytopes possess a regular unimodular triangulation. 
\end{abstract}

\section*{Introduction} 

Let $G$ be a finite connected simple graph and $P_G$ be the edge polytope of $G$. 
The combinatorial structure of $P_G$, especially which type of triangulations $P_G$ admits, is an interesting problem, and many research studies have been done on this topic (see \cite[Chapter 5]{dojo} and references therein). 
In \cite{Ohsugi}, Ohsugi obtained a necessary and sufficient condition for $P_G$ to possess a regular unimodular triangulation 
(there exists a monomial order such that the initial ideal of the toric ideal of the graph $G$ is generated by squarefree monomials).  
However, this condition is not so easy to apply to a given graph just by looking at the graph.  

In this paper, for a graph $G$, we will give several criteria for the existence of a regular unimodular triangulation of $P_G$ in terms of some simple data of the graph. 
We also apply our criteria to some examples and show that their edge polytopes possess a regular unimodular triangulation. 

The contents of this paper are as follows. 
In Section 1, we review the definitions of and some basic results on the graphs in \cite{Ohsugi}. 
In Section 2, we give several slightly different criteria for $P_G$ to possess a regular unimodular triangulation.  
In Section 3, we show some examples to which our criteria are applicable.

\section{Preliminaries} 
A matrix $A=(a_{ij})_{1\le i \le d, 1\le j \le n} \in {\bf Z}^{d\times n}$ is called a {\it configuration matrix} if there exists {\bf c} $\in {\bf R}^d$ such that
$${\bf a}_j\cdot{\bf c}=1, \ 1\le j\le n$$
where ${\bf a}_j$ is a column vector of $A$.

Let $\mathcal{A}=[{\bf a}_1,\dots ,{\bf a}_n]\in{\bf Z}^{d\times n}$ be a configuration matrix.
Let $\Delta$ be a collection of simplices whose vertices belong to a configuration matrix $\mathcal{A}$.
Then, $\Delta$ called a {\it covering} of $\mathcal{A}$ if
$${\rm CONV}(\mathcal{A})=\bigcup_{F\in\Delta}F$$
holds.
In addition, if a covering $\Delta$ of a configuration matrix $\mathcal{A}$ is a simplicial complex, then it is called a {\it triangulation} of $\mathcal{A}$.
For a configuration matrix $\mathcal{A}=[{\bf a}_1,\dots ,{\bf a}_n]\in{\bf Z}^{d\times n}$, let
$${\bf Z}\mathcal{A}=\left\{\sum_{i=1}^n z_i{\bf a}_i \ : \ z_i\in{\bf Z}\right\}\subset{\bf Z}^d.$$

Let $\mathcal{B}\subset\{{\bf a}_1,\dots ,{\bf a}_n\}$ be the vertex set of a maximal simplex $\sigma \in \Delta$ in a covering (triangulation) $\Delta$ of $\mathcal{A}$.
Suppose that the rank of a configuration matrix $\mathcal{A}\in{\bf Z}^{d\times n}$ is equal $d$.
Let $\delta$ be the greatest common divisor of all $d \times d$ minors of $\mathcal{A}$.
Then, the {\it normalized volume} of $\sigma$ is defined by
$${\rm VOL}(\sigma)=\dfrac{|\det (\mathcal{B})|}{\delta}.$$

A covering (triangulation) $\Delta$ of $\mathcal{A}$ is said to be {\it unimodular} if the normalized volume of any maximal simplex in $\Delta$ is equal to 1.
For a configuration matrix $\mathcal{A}=[{\bf a}_1,\dots ,{\bf a}_n]\in{\bf Z}^{d\times n}$ and a vector ${\bf w}=[w_1,\dots,w_n]\in{\bf Q}^n$, let $\Delta_{\bf w}$ be the set of all convex polytopes ${\rm CONV}(\{{\bf a}_{i_1},\dots {\bf a}_{i_r}\})$ satisfying the following condition:
$${\rm There \ exists} \ {\bf c} \in {\bf Q}^d {\rm such \ that} \begin{cases} {\bf a}_j\cdot{\bf c}=w_j, \ \ j\in\{i_1,\dots,i_r\}, \\ {\bf a}_j\cdot{\bf c}<w_j, \ \ j\notin\{i_1,\dots,i_r\}. \end{cases}$$

A triangulation $\Delta$ of a configuration matrix $\mathcal{A}$ is said  to be {\it regular} if there exists ${\bf w}\in {\bf Q}^d$ such that $\Delta=\Delta_{\bf w}$.

Let $t_1,t_2,\dots, t_d$ be variables.
Let $\mathcal{A}=(a_{ij})_{1\le i \le d, 1\le j \le n} \in {\bf Z}^{d\times n}$ be a configuration matrix.
To each column vector
$${\bf a}_j=\begin{pmatrix} a_{1_j} \\ a_{2_j} \\ \vdots \\ a_{d_j} \end{pmatrix},$$
we associate the monomial
$${\bf t}^{{\bf a}_j}=t_1^{a_{1_j}}t_2^{a_{2_j}}\dots t_d^{a_{d_j}}$$
with allowing negative powers.
Let $K$ be a field and let $K[{\bf x}]=K[x_1,x_2,\dots,x_n]$ be a polynomial ring in $n$ variables over $K$.
If $f=f(x_1, x_2, \dots, x_n)\in K[{\bf x}]$, then we define $\pi(f)$ by setting
$$\pi(f)=f({\bf t}^{{\bf a}_1},{\bf t}^{{\bf a}_2},\dots,{\bf t}^{{\bf a}_n}).$$
Let
$$K[\mathcal{A}]=\{\pi(f) : f\in K[{\bf x}]\}.$$
We say that $K[\mathcal{A}]$ is the {\it toric ring} of $\mathcal{A}$.
In general, a configuration matrix $\mathcal{A}$ satisfies ${\bf Z}_{\ge 0}\mathcal{A}\subset{\bf Z}\mathcal{A}\cap{\bf Q}_{\ge 0}\mathcal{A}$.
The toric ring $K[\mathcal{A}]$ is said to be {\it normal} if it satisfies ${\bf Z}_{\ge 0}\mathcal{A}={\bf Z}\mathcal{A}\cap{\bf Q}_{\ge 0}\mathcal{A}$.

With respect to the normality of the toric ring $K[\mathcal{A}]$, the existence of unimodular triangulations and uniodular coverings of $\mathcal{A}$ plays an important role.

\begin{prop}[{\cite[Theorem 5.6.7]{dojo}}]
If a configuration matrix $\mathcal{A}$ has a unimodular covering, then the tric ring $K[\mathcal{A}]$ is normal.
\end{prop}

Let $G=(V,E)$ be a finite connected simple graph, where $V=\{1,2, \dots, d \}$ is the vertex set and $E = \{ e_1, \dots,e_n \}$ is the set of edges. 
Here, a graph is called {\it simple} if it has no loop and no multiple edges. 
For each edge $e = \{i, j \} \in E$, we set $\rho(e) := {\bf e}_i + {\bf e}_j \in {\bf Z}^d$, where ${\bf e}_i$ is the $i$th unit coordinate vector in ${\bf R}^d$. 
We call the convex hull $P_G \subset {\bf R}^d$ of the finite set $\{ \rho(e) \,| \, e \in E \}$ as the {\it edge polytope} of $G$.     

A sequence $\Gamma=(e_{j_1},\dots ,e_{j_r})$ of edges of a finite graph $G$ is called a {\it walk} of length {\it r} if $\Gamma$ satisfies
$$e_{j_1}=\{i_1,i_2\}, e_{j_2}=\{i_2,i_3\}, \dots ,e_{j_r}=\{i_r,i_{r+1}\}.$$
In addition, if $i_{r+1}=i_1$, then $\Gamma$ is called a {\it closed walk} of length {\it r}.
A closed walk of even length is called an {\it even closed walk}.
If $i_{r+1}=i_{1}$ and $i_1,\dots ,i_r \ (r \ge 3)$ are distinct, then $\Gamma$ is called a {\it cycle} of length {\it r}.
An edge that joins two vertices of a cycle but is not itself an edge of the cycle is a {\it chord} of that cycle (\cite{G}).
A cycle $C$ in a graph is called {\it minimal} if $C$ possesses no chord.

Let $C$ be an odd cycle contained in $G$. 
Let $c$ be a chord of $C$. 
Then, $c$ divides $C$ into two cycles, where one is an odd cycle and the other is an even cycle. 
We call the even cycle {\it the even closed walk of the chord $c$ in $C$}. 
In the even closed walk $\Gamma$ of the chord $c$ in $C$, we require that $c$ be an even-numbered edge of $\Gamma$. 

Let $G$ be a finite connected simple graph on the set of vertex $\{1,\dots,d\}$.
Let $K[t_1,\dots, t_d]$ denote the polynomial ring in $d$ indeterminates over a field $K$ and let $K[G]$ be the subalgebra of $K[t_1,\dots,t_d]$ generated by all quadratic monomials $t_it_j$ such that $\{i,j\}$ is an edge of $G$.
The affine semigroup ring $K[G]$ is called the {\it edge ring} of $G$. 

Let $(C_1, C_2)$ be a pair of disjoint odd cycles in $G$ (namely, the odd cycles $C_1$ and $C_2$ have no common vertex) and $b$ be a bridge of this pair.  
Here, a bridge $b$ of the pair $(C_1, C_2)$ is an edge $b = \{ i, j \}$, where $i$ is a vertex of $C_1$ and $j$ is a vertex of $C_2$ or vice versa.   
Then, {\it the even closed walk of $b$ in $(C_1, C_2)$} is the closed walk $(C_1, b, C_2, -b)$. 
In this notation, $-b$ means the oppositely directed edge of $b$ and the cycle $C_1$ starts from the vertex $C_1 \cap b$ and ends at the same vertex. 
The same holds for $C_2$.  
We note that, in the even closed walk $\Gamma$ of the bridge $b$ in $(C_1, C_2)$, $b$ appears twice as an even-numbered edge of $\Gamma$. 

A {\it Fulkerson–Hoffman–McAndrew (FHM) graph} (\cite{FHM}) is a finite connected simple graph such that any pair of disjoint odd cycles has a bridge. 
A {\it fundamental FHM graph} (\cite{FHM}) is an FHM graph that has at least one pair of disjoint odd cycles. 

It is also known that the normality of edge polytopes is characterized by the following condition.
\begin{prop}[{\cite[Corollary 2.3]{OH}}]
Let $G$ be a finite connected simple graph.
Then the following conditions are equivalent:
\begin{enumerate}
\item[{\rm (i)}] the edge ring $K[G]$ is normal;
\item[{\rm (ii)}] the edge polytope $P_G$ possesses a unimodular covering;
\item[{\rm (iii)}] the graph $G$ is an FHM graph.
\end{enumerate}
\end{prop}

The following is a basic fact about the fundamental FHM graph (\cite[Corollary 2.3]{OH}, \cite[Proposition 3.4]{Ohsugi}, \cite{OH2}, and \cite{A}).

\begin{prop} Let $G$ be a finite connected simple graph. 
\begin{enumerate}
\item[{\rm (i)}] If the edge polytope $P_G$ possesses a regular unimodular triangulation, then $G$ is an FHM graph.  
\item[{\rm (ii)}] If $G$ possesses no pair of disjoint odd cycles, then $P_G$ possesses a regular unimodular triangulation. 
\item[{\rm (iii)}] There exists an example of an edge polytope $P_G$ of a fundamental FHM graph $G$ that possesses no regular unimodular triangulation.
\end{enumerate}
\end{prop}

Thus, we focus on the fundamental FHM graph hereafter. 
We will review the necessary and sufficient condition for $P_G$ to have a regular unimodular triangulation.  

Let $G$ be a fundamental FHM graph. 
Suppose $G$ possesses $p$ pairs of disjoint odd cycles $\Pi_1= (C_1,C'_1), \dots, \Pi_p= (C_p,C_p')$. 
For each $i$ ($1 \leq i \leq p$), let $\{ b^i_j \mid 1 \leq j \leq q_i \}$ be the set of bridges of $\Pi_i$ and the chords of $C_i$ or $C'_i$.
Let $\Gamma^i_j = (e_{i_1}e_{i_2} \dots e_{i_{2r}})$ be the even closed walk of $b^i_j$, where the bridge or chord is even numbered. 

Now, we define the open half-space $H_{b^i_j}$ by
\begin{equation}
H_{b^i_j}:= \{ (x_1, \dots,x_n)\in {\bf R}^n \mid \sum_{k=1}^r x_{i_{2k-1}} > \sum_{k=1}^r x_{i_{2k}}\}. 
\end{equation}

Furthermore, we set $ W:= \bigcap_{i=1}^p(\bigcup_{j=1}^{q_i} H_{b^i_j})$.
The following result is our starting point. 

\begin{prop}[{\cite[Theorem 3.5]{Ohsugi}}]
The edge polytope $P_G$ possesses a regular unimodular triangulation if and only if $W \ne \phi$.  
\end{prop}

\section{Criteria for the existence of a regular unimodular triangulation}

Let $G$ be a fundamental FHM graph. 
In this section, we will give four criteria for the edge polytope $P_G$ to possess a regular unimodular triangulation in terms of the simple data of the graph $G$. 
Our criteria are based on the existence of special bridges in each pair of disjoint odd cycles. 
Let $\Pi_1,\dots, \Pi_p$ be all the pairs of disjoint odd cycles in $G$ as before and $\{b^1,\dots, b^p\}$ be the set of bridges, where $b^i$ is the bridge of $\Pi_i$. 
Let $\Gamma_i:=(e_{i_1},\dots, e_{i_{2s+1}}, b^i, e_{j_1},\dots, e_{j_{2t+1}}, -b^i).$

Now, we define 
$$\alpha_i:=|\{b^1,\dots, b^p\}\cap\{e_{i_2},e_{i_4},\dots, e_{i_{2s}},e_{j_2},e_{j_4},\dots,e_{j_{2t}}\}|,$$
$$\beta_i:=|\{b^1,\dots, b^p\}\cap\{e_{i_1},e_{i_3},\dots, e_{i_{2s+1}},e_{j_1},e_{j_3},\dots,e_{j_{2t+1}}\}|.$$ 
Furthermore, we set $a_i:=2+\alpha_i-\beta_i$.

\begin{thm}
Work with the same notation as above.
The edge polytope of a fundamental FHM graph $G$ possesses a regular unimodular triangulation if 
it has a set of bridges $\{b^1, \dots,b^p \}$ ($b^i$ is the bridge of $\Pi_i$) that satisfies the following condition:
For each $i$, $a_i \geq 0$ holds, and furthermore, the number of $\Gamma_i's$ such that $a_i = 0$ is at most 2.
\end{thm}

\noindent
{\it Proof.} We first rewrite $W$ in Proposition 1.4 by the distributive law as follows:  
$$W =  \bigcap_{i=1}^p(\bigcup_{j=1}^{q_i} H_{b^i_j}) = \bigcup_{j_1, \dots,j_p} (H_{b^{1}_{j_1}} \cap \dots \cap H_{b^{p}_{j_p}}),$$
where $j_k$ satisfies $1 \leq j_k \leq q_k$.
We set
$$C_b= C_{\{b^1_{j_1}, \dots, b^p_{j_p}\}}:= H_{b^{1}_{j_1}} \cap \dots \cap H_{b^{p}_{j_p}}$$ and call $C_b$ as the open cone of $b= \{b^1_{j_1}, \dots, b^p_{j_p}\}$.
Thus, $W \ne \phi$ is equivalent to that there is a set of bridges $b = \{b^1, \dots,b^p\}$ ($b^i$ is a bridge of $\Pi_i$) such that $C_b$ is non-empty.  

For each $i$, let $\Gamma_i$ be the even closed walk of $b^i$ and $f_i>0 $ be the inequality (1) defined by $b^i$. 
We denote by the same $f_i$ an $n$-dimension vector that consists of the coefficients of the left-hand side (LHS) of the inequality $f_i>0$.  
We note if the bridge $b^i$ is equal to an edge $e_j$,  if the $j$th component $f_i[j]$ of the vector $f_i$ is $-2$, and if the other edge $e_k$ is contained in $\Gamma_i$, $f_i[k]= +1$ (respectively $-1$) if $e_k$ is an odd (respectively even)-numbered edge of $\Gamma_i$.  
The other components of $f_i$ are 0.

We define the standard weight vector $w \in {\bf R}^n$ of $C_b$ as follows. 
If there exists $i$ such that $f_i[k]= -2$, then we set $w[k]:= -1$. 
The other components of $w$ are 0.
We note if $a_i$ is equal to $(f_i, w)$ (inner product) for each $i$.
\begin{enumerate}
\item[{\rm (i)}] Suppose $a_i>0$ for any $i$.
Since $(f_i, w)>0$ for any $i$, then $w \in W$

\item[{\rm (ii)}] Suppose $a_j=0$ and $a_i>0 \ (i\neq j)$.
Let $b^j$ be a bridge of $\Gamma_j$ and $b^j = e_l$.
Let $w':=w+(-1/10\cdot{\bf e}_l$), where ${\bf e}_l$ is a unit vector.
Now, we consider $(f_j,w')=(f_j,w)+(f_j, -1/10\cdot{\bf e}_l)$.
By the assumption, $(f_j, w)=a_j=0$.
Moreover, we obtain $(f_j, -1/10\cdot{\bf e}_l)=1/5$.
Therefore, $(f_j, w')=(f_j,w)+(f_j, -1/10\cdot{\bf e}_l)=1/5>0$.
On the other hand, let $b^k$ be a bridge of $\Gamma_i$ and $b^k = e_m$.
Let $w':=w+(-1/10\cdot{\bf e}_m)$.
Next, we consider $(f_i,w')=(f_i,w)+(f_i, -1/10\cdot{\bf e}_m)$.
By the assumption, $(f_i,w)=a_i>0$.
Moreover, we obtain $(f_i, -1/10\cdot{\bf e}_m)=1/5$.
Therefore, $(f_i,w')=(f_i,w)+(f_i, -1/10\cdot{\bf e}_m)>0$.

\item[{\rm (iii)}] Suppose $a_j=a_k=0$ and $a_i>0$ \ $(i\neq j, i\neq k)$.
There exists at least an edge $e_l$ in $\Gamma_j$ that is not contained in $\Gamma_k$.
On the other hand, there exists at least an edge $e_m$ in $\Gamma_k$ that is not contained in $\Gamma_j$.
Let $v$ be a vector that satisfies the following condition:
$v[l]=1/10$ (respectively $-1/10$) if $e_l$ is odd numbered (respectively even numbered) in $\Gamma_j$.
$v[m]=1/10$ (respectively $-1/10$) if $e_m$ is odd numbered (respectively even numbered) in $\Gamma_k$.
The other components of $v$ are 0.
Let $w':=w+v$.
Now, we consider $(f_i, w')=(f_i,w)+(f_i,v)$.
By the assumption, $(f_i,w)=a_i>0$.
On the other hand, we obtain $(f_i,v)\ge -3/10$.
Here, since $a_i\in{\bf Z}_{>0}$, then $(f_i,w)=a_i\ge 1$.
Therefore, $(f_i, w')=(f_i,w)+(f_i,v)\ge 7/10>0$.
Next, we consider $(f_j,w')=(f_j,w)+(f_j, v)$.
By the assumption, $(f_j, w)=a_j=0$.
Moreover, we obtain $(f_j, v)= 1/10$ or $1/5$.
Therefore, $(f_j, w')=(f_j,w)+(f_j, v)>0$.
Work with the same discussion as above.
We obtain $(f_k, w')=(f_k,w)+(f_k, v)>0$.
\end{enumerate}

We have the following corollaries.

\begin{cor}
Work with the same notation as above.
The edge polytope of a fundamental FHM graph $G$ possesses a regular unimodular triangulation if it has a set of bridges $\{b^1, \dots,b^p \}$ ($b^i$ is the bridge of $\Pi_i$) that satisfies the following condition: For each $i$, $a_i > 0$ holds. 
\end{cor}

\begin{cor}
The edge polytope of a fundamental FHM graph $G$ possesses a regular unimodular triangulation if it has a set of bridges $\{b^1, \dots,b^p \}$ ($b^i$ is the bridge of $\Pi_i$) that satisfies the following condition: for each even closed walk $\Gamma_i$ of $b^i$, the number of the other bridges $b^j$ contained in $\Gamma_i$ is at most 2, and, furthermore, the number of $\Gamma_i's$ that contain exactly two other bridges is at most 2. 
\end{cor}

\begin{cor}  
The edge polytope of a fundamental FHM graph $G$ possesses a regular unimodular triangulation if it has a set of bridges $\{b^1, \dots,b^p \}$ ($b^i$ is the bridge of $\Pi_i$) that satisfies the following  condition: each even closed walk of the bridge $b^i$ contains at most one other bridge $b^j$.  
\end{cor}

We note that the narrowest condition is Corollary 2.4, whereas the broadest is Theorem 2.1.
However, Corollary 2.4 is the easiest to check graphically.

\begin{chu}
\begin{enumerate}
\item[{\rm (i)}] {\rm In Theorem 2.1, if there exist more than two $i's$ such that $a_i=0$, the following result holds. Suppose $a_i = 0$ for $i= i_1, \dots,i_r$ ($r \geq 3$) and $a_i>0$ for the other $i's$. 
Let $H \subset {\bf R}^n$ be the hyperplane defined by $\sum_{j=1}^n w[j]x_j=0$. 
If the convex cone generated by $f_{i_1}, \dots,f_{i_r}$ in $H$ is strongly convex i.e., for a convex cone $P$, $P\cap-P=\{{\bf 0}\}$, $W$ is not empty.
The proof is the same as that of Theorem 2.1. 
Namely, thanks to this condition, we can vary $w$ slightly to get a new weight $w'$ such that $(f_i,w')>0$ for any $i$. 
However, this condition is not clear at all just by looking at the graph. }

\item[{\rm (ii)}] {\rm More generally, let $C(f_1,\dots, f_p)$ be an open cone in ${\bf R}^n$ defined by $p$ linear homogeneous inequalities $f_i > 0$ ($1 \leq i \leq p$). 
Then, $C(f_1,\dots, f_p) \ne \phi$ holds if and only if the dual cone $C(f_1,\dots, f_p)^{\vee}= {\bf R}_{\geq 0}f_1+ \dots +{\bf R}_{\geq 0} f_p$ of $C(f_1,\dots, f_p)$ is strongly convex ($f_i$ is the coefficient vector of the LHS of the inequality). 
It is difficult to determine whether $C(f_1,\dots f_p)^{\vee}$ is strongly convex or not just by looking at the graph.}
\item[{\rm (iii)}] {\rm The edge polytope of the following graph does not possess the regular unimodular triangulations (Example 3.2 in \cite{01}). 
Moreover, there exist three $i's$ such that $a_i=0$. 
Therefore, we cannot improve the condition of Theorem 2.1 such that ``the number of $\Gamma_i's$ such that $a_i = 0$ is at most 3'' }

\begin{center}
\includegraphics[width=3.5cm]{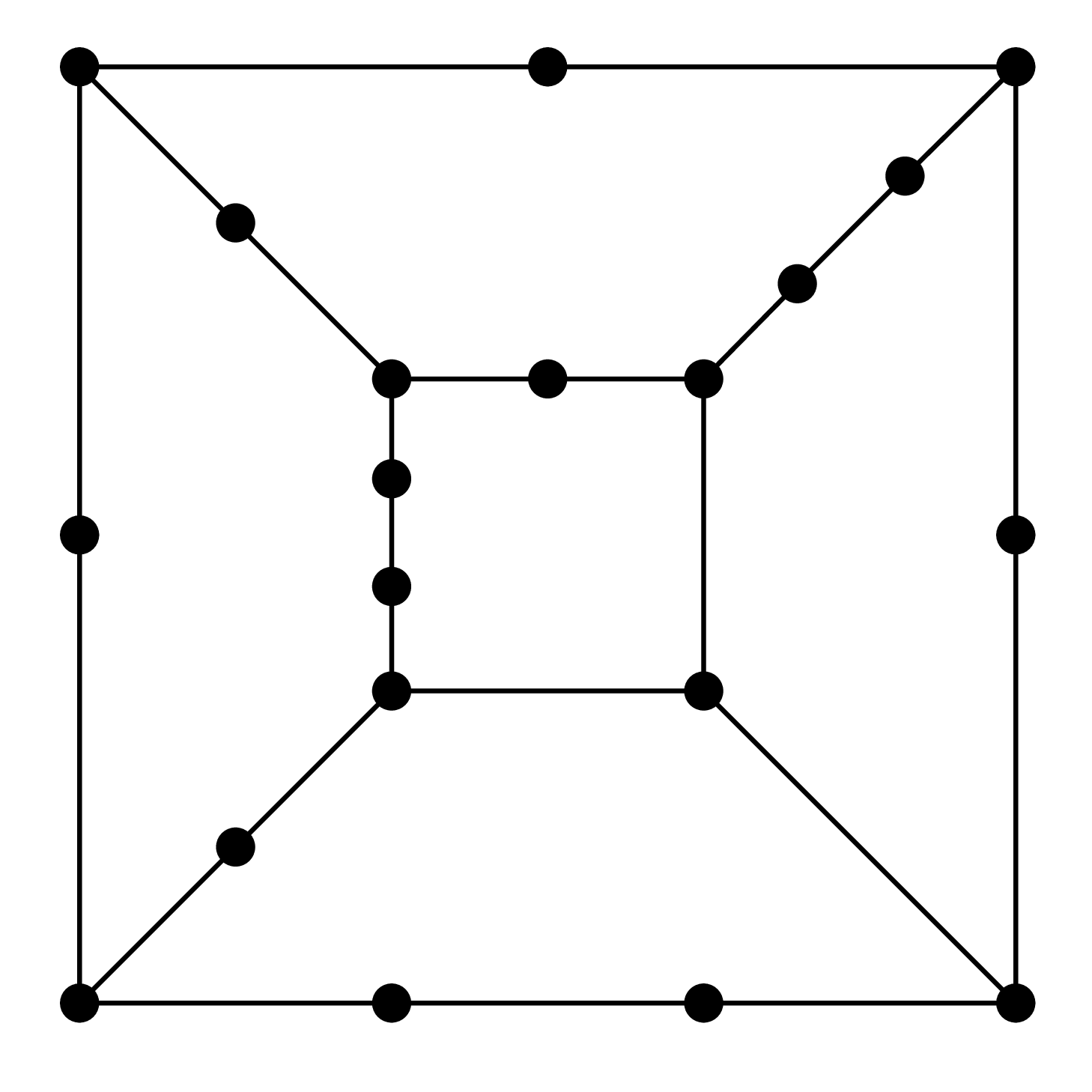}
\end{center}

\end{enumerate}
\end{chu}  

\section{Applications}

We first apply our criteria to the complete graph $G = K_6$ with six vertices. 
It is known that $P_{K_d}$ possesses a regular unimodular triangulation for any $d$ (see \cite{B}).  
Moreover, it is known that an edge polytope of a gap-free graph or a complete multipartite graph possesses a regular unimodular triangulation too (see \cite{D} and \cite{OH3}).

\begin{chu}
{\rm The complete graph $K_6$ satisfies the condition of Corollary 2.2, but does not satisfy the condition of Corollary 2.3.}
\end{chu}

We finally show several other examples that satisfy our criteria.

\begin{rei}
{\rm The following five types of graphs satisfy the condition of Corollary 2.4. 
More precisely, in the graphs $A_{m,n}, B_{m,n}$, and $C_{m_1,m_2,n_1,n_2}$, all the pairs of disjoint odd cycles (triangles) have a bridge $b$ in common, and, thus, there are no other bridges contained in the even closed walk of $b$. 

$D_{m_1,m_2,m_3,m_4}$ has a set of bridges $\{b_1,b_2\}$ where any disjoint pair has a bridge in this set, and the even closed walk of $b_i$ ($i=1,2$) contains (exactly) one other bridge.  
$E_{m_1,m_2,m_3}$ has a set of three bridges $\{ b_1,b_2,b_3 \}$ where any disjoint pair has a bridge in this set, and there are no other bridges contained in the even closed walk of $b_i$ ($i=1,2,3$).  } 
\end{rei}

\noindent
\includegraphics[width=7.5cm]{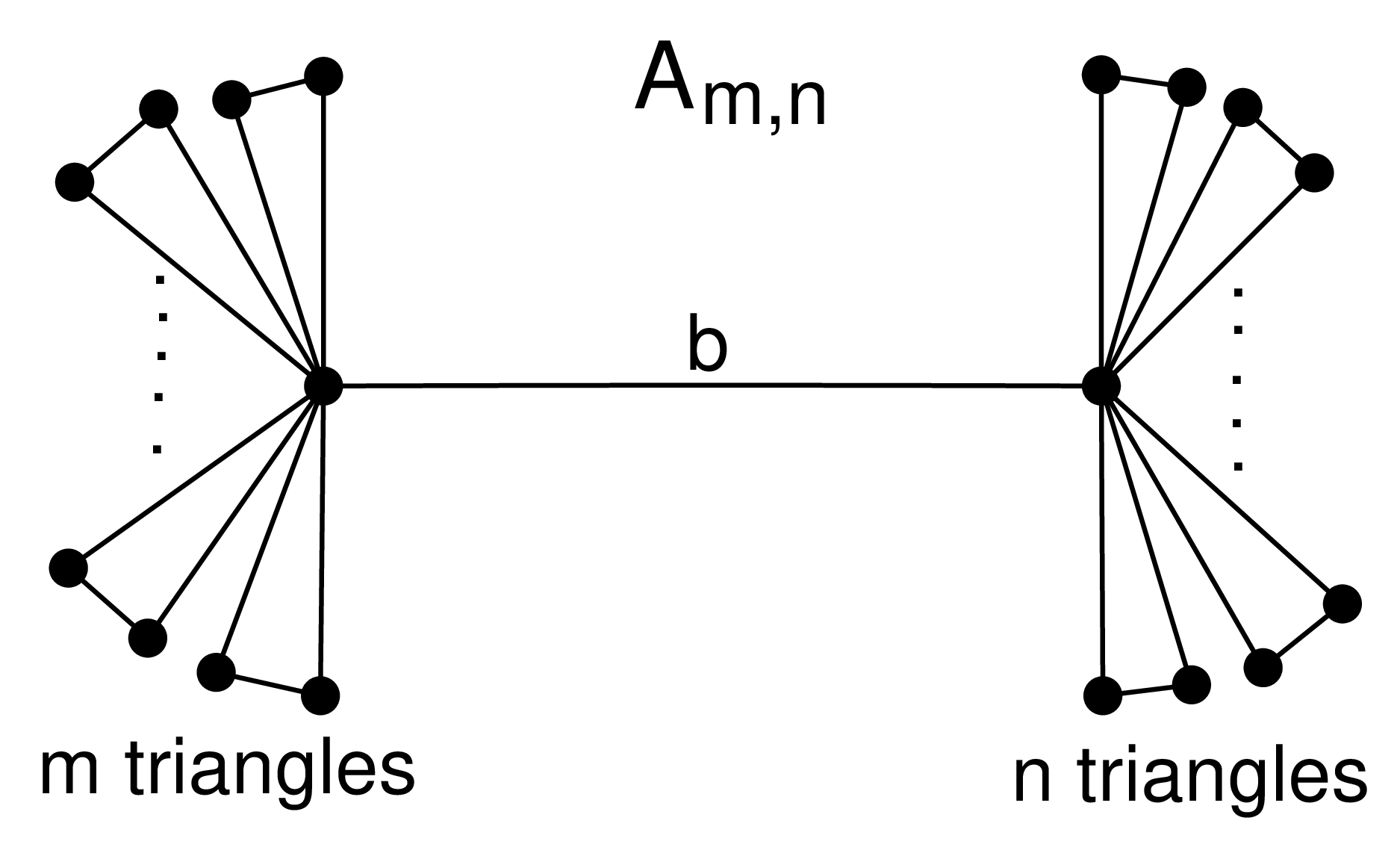} \ \ \ \ \includegraphics[width=7.5cm]{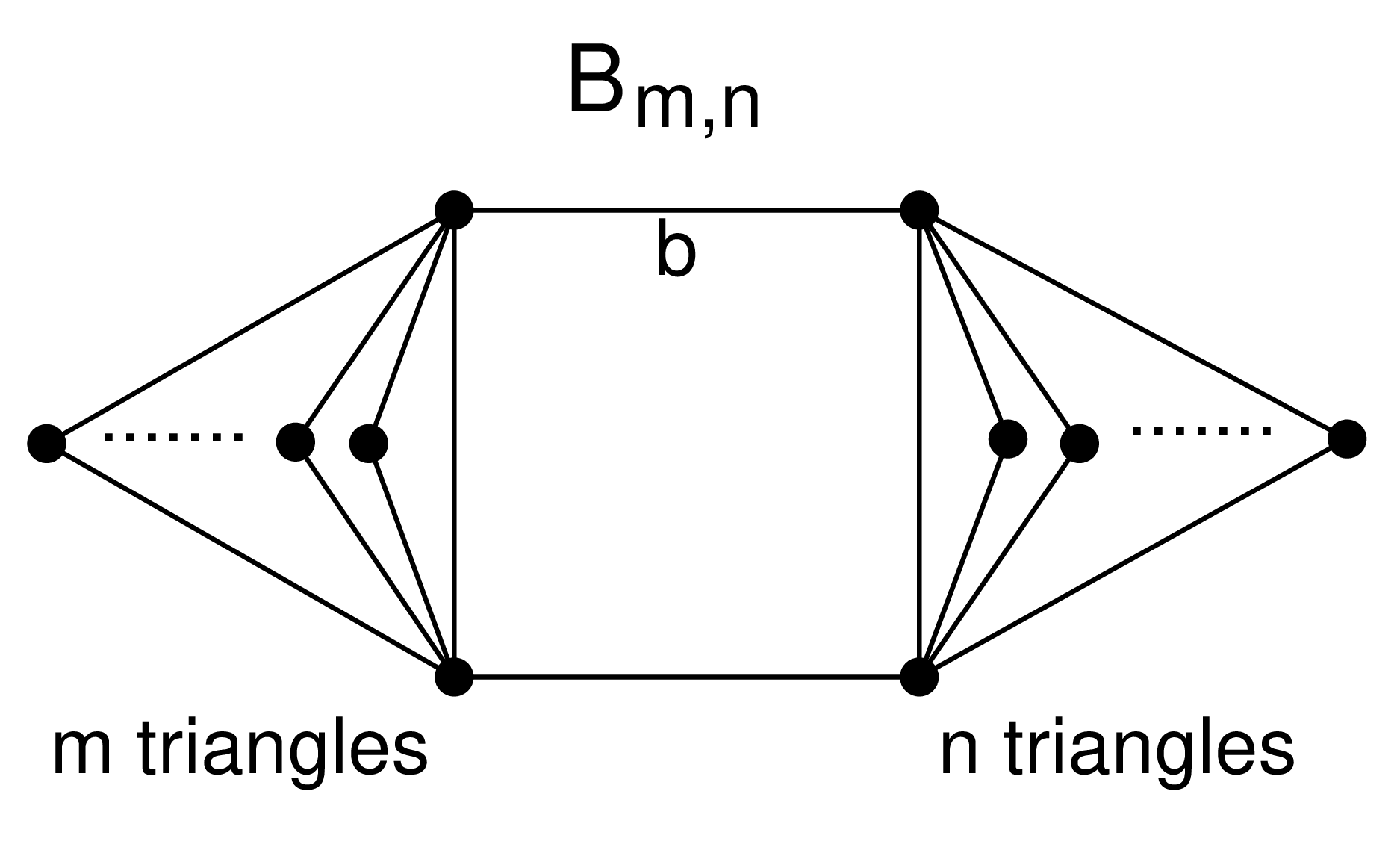}\\\\
\includegraphics[width=7.5cm]{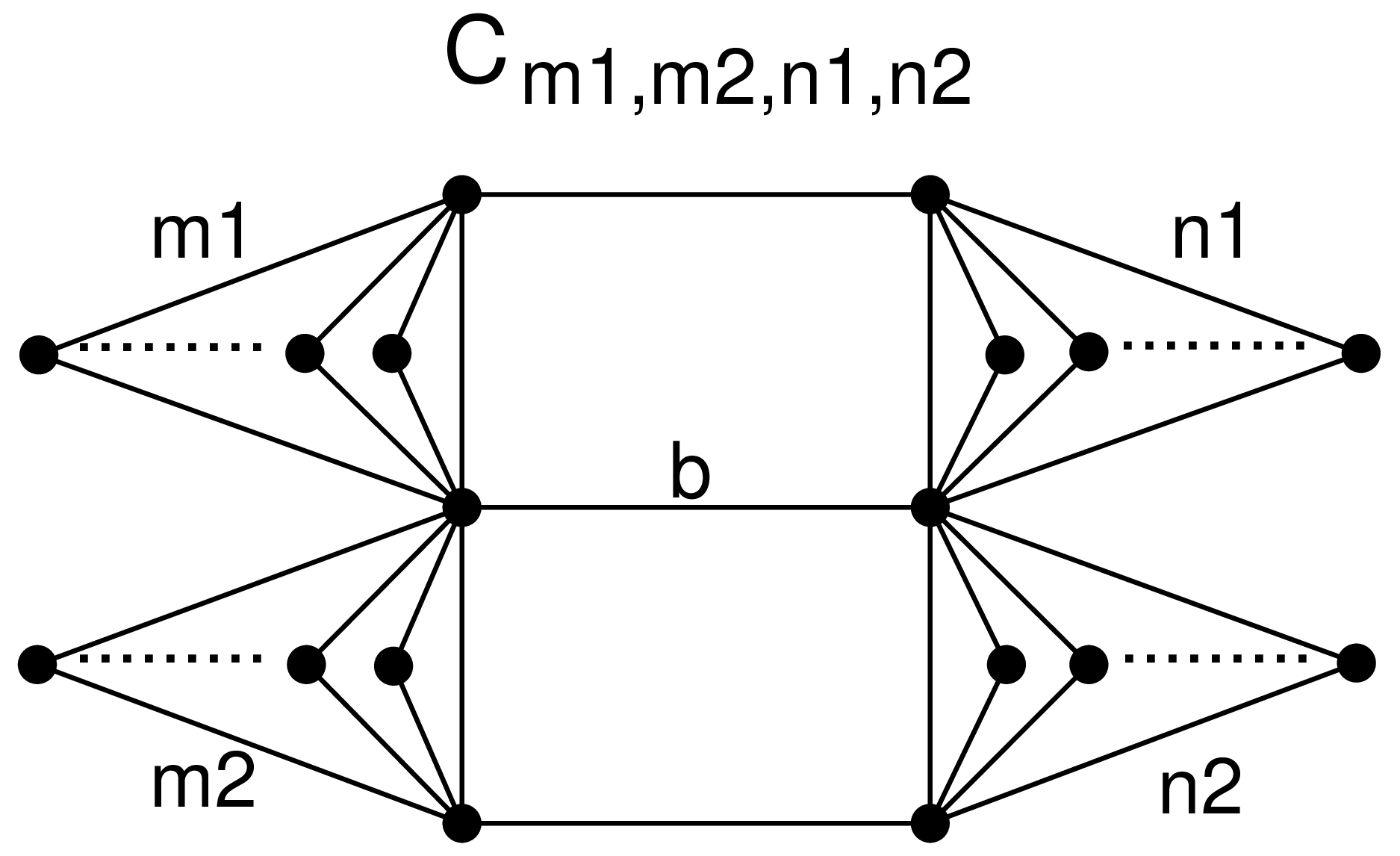} \ \ \ \ \includegraphics[width=8cm]{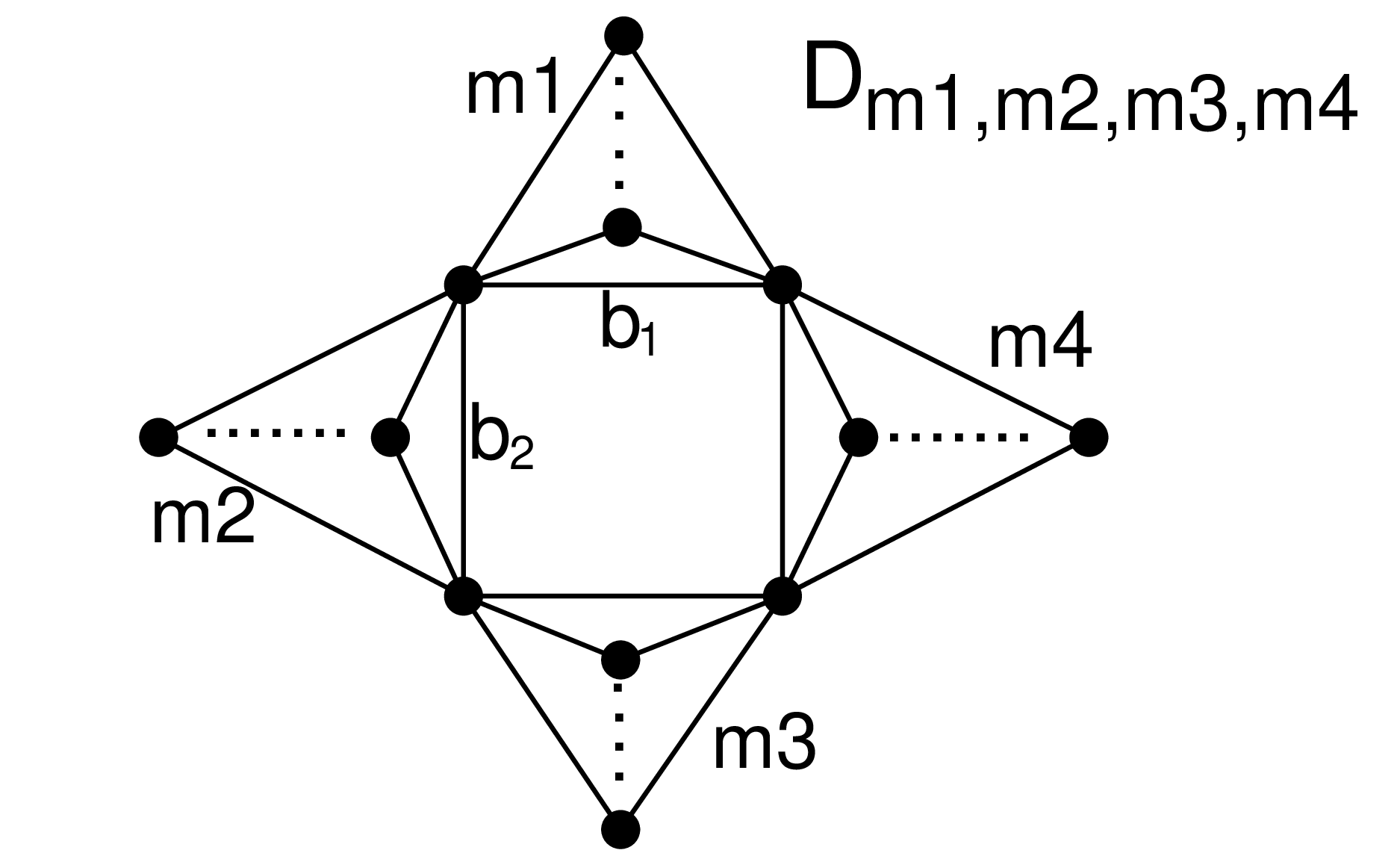}\\\\
\includegraphics[width=8cm]{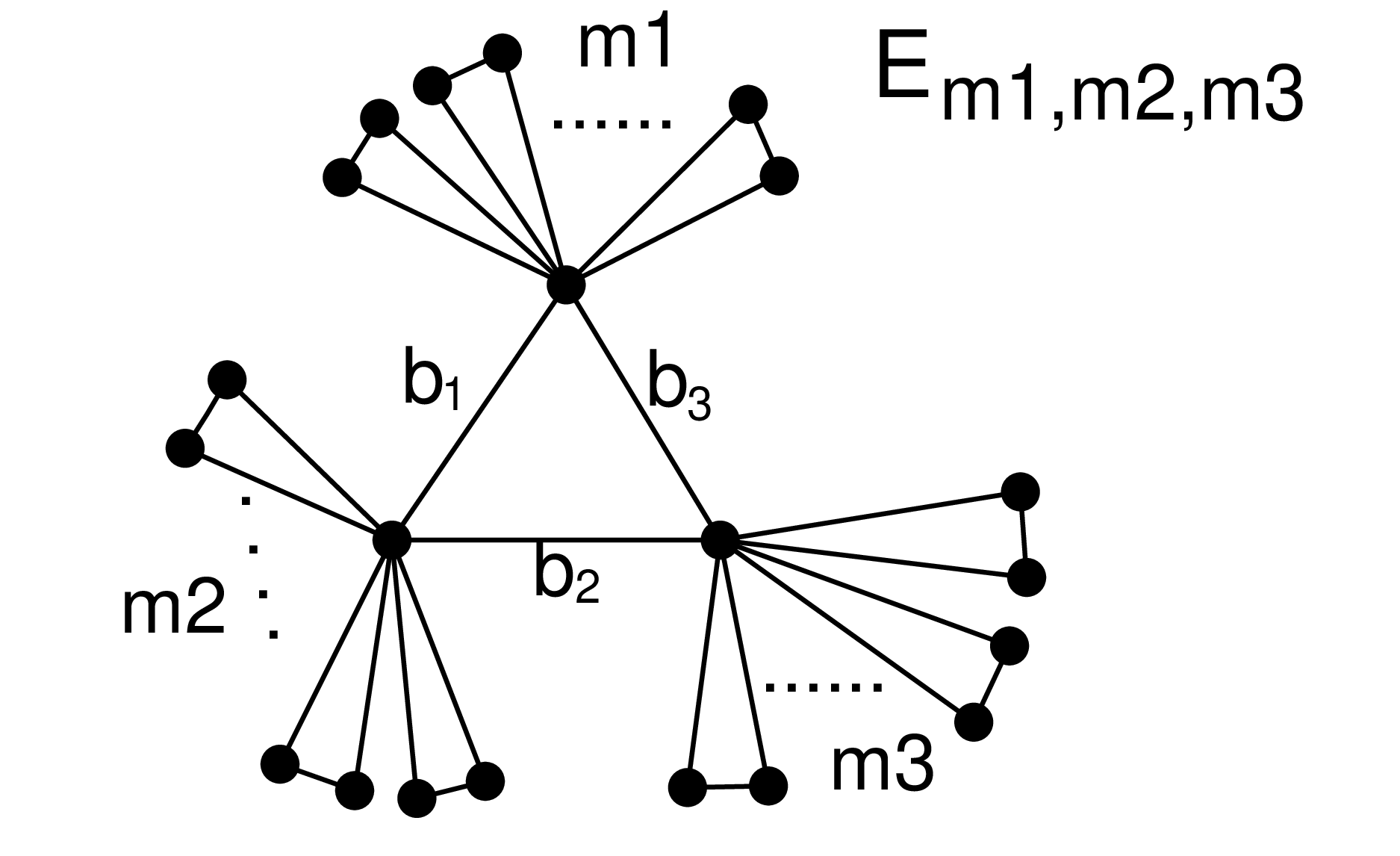}

\begin{rei}
{\rm The following two types of graphs satisfy the condition of Corollary 2.2, but not that of Corollary 2.3. 
$F_{m_1,m_2,m_3,m_4}$ has a minimal set of six bridges $\{ b_i \,| \, 1 \leq i \leq 6 \}$ where any disjoint pair has a bridge in this set, and $G_{m_1,m_2,m_3,m_4,m_5}$ has a minimal set of ten bridges $\{ b_i \,| \, 1 \leq i \leq 10 \}$.}   
\end{rei}
\noindent
\includegraphics[width=8cm]{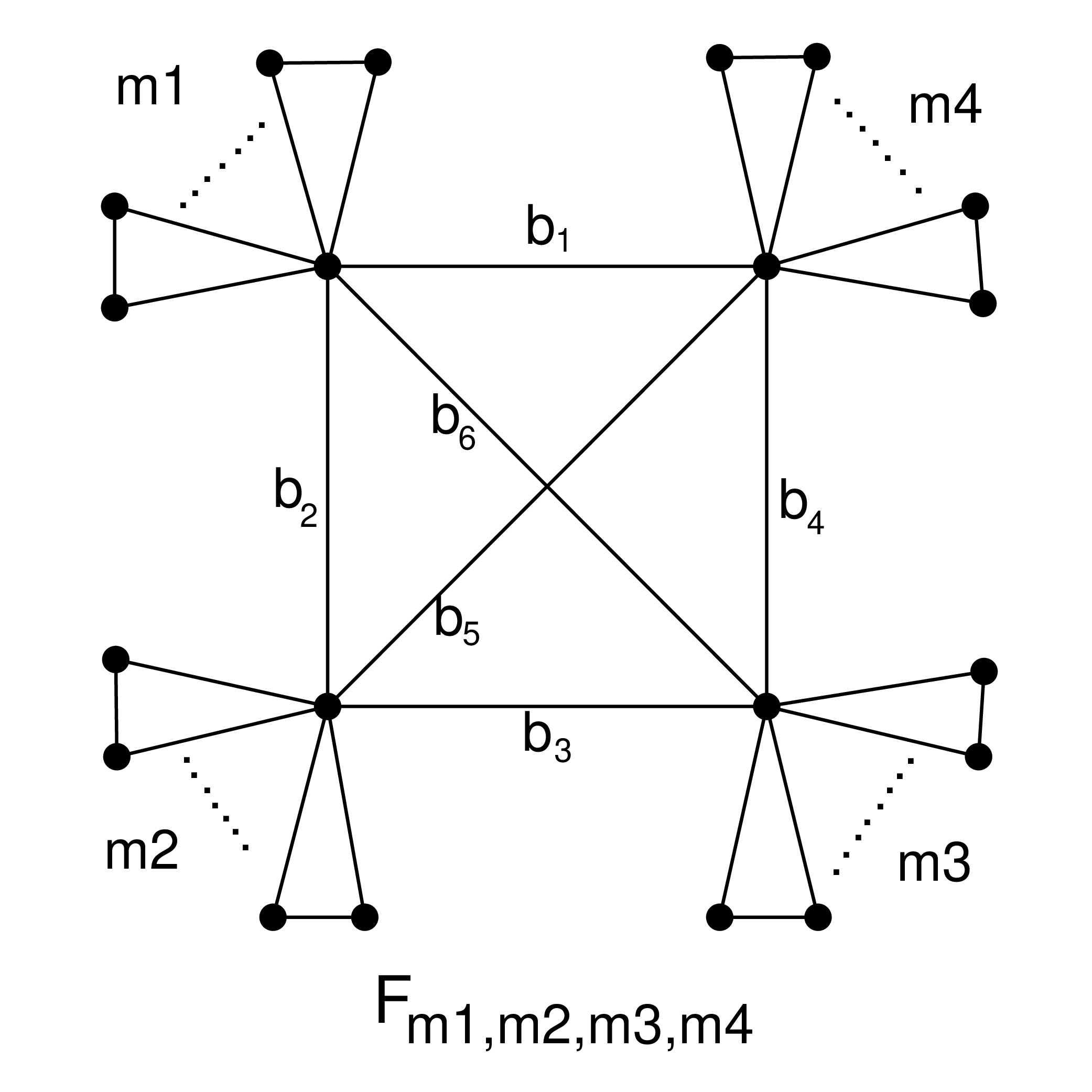}\includegraphics[width=8cm]{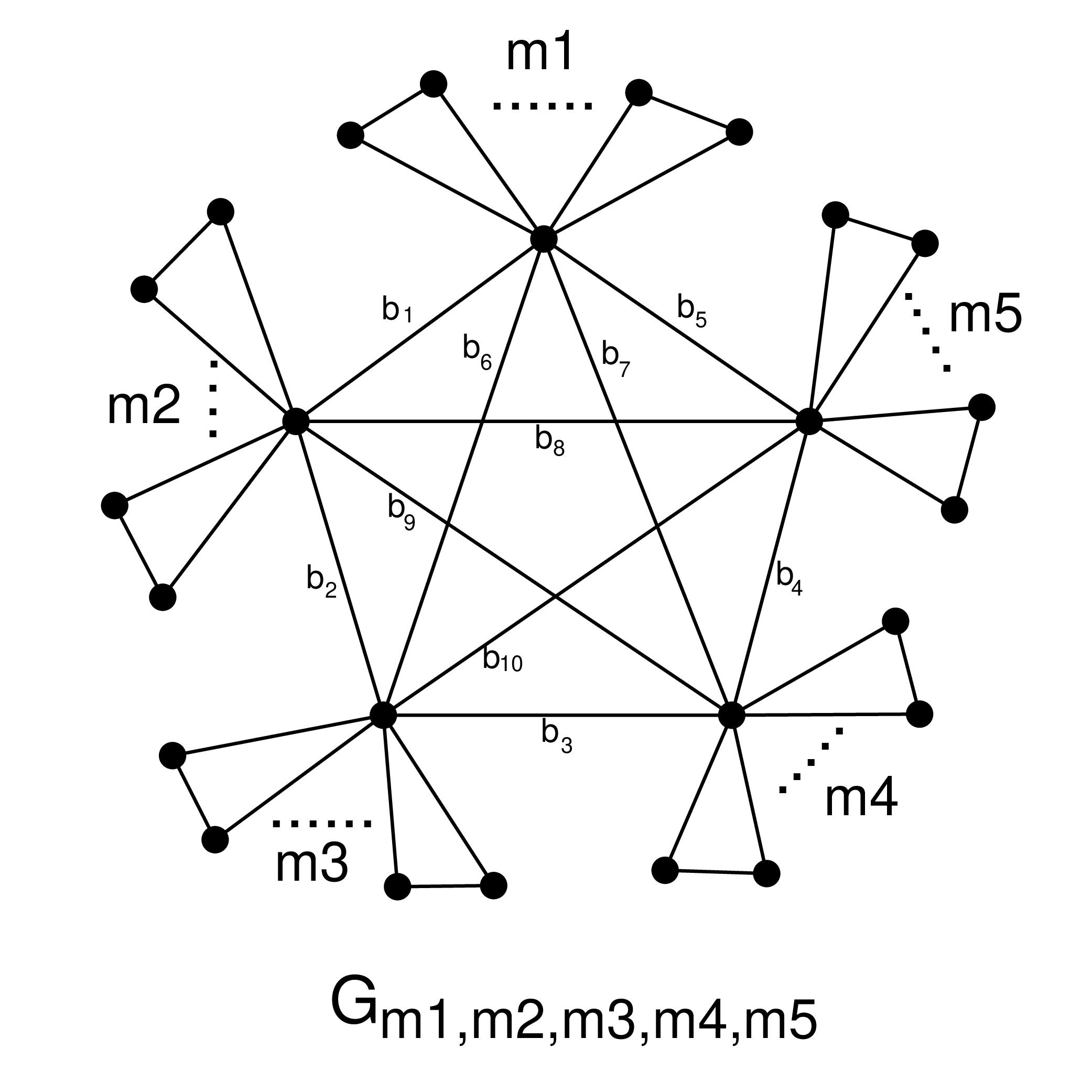}

\begin{rei}
{\rm The following graph satisfies the condition of Theorem 2.1.
Moreover, there exist just two $i's$ such that $a_i=0$.
The following graph has a minimal set of seven bridges $\{ b_i \,| \, 1 \leq i \leq 7 \}$ where any disjoint pair has a bridge in this set.
When $\Gamma_1=(e_{24},e_{16},e_{15},b_7,b_1,b_5,b_6,e_{22},b_2,-b_7)$, then $\alpha_1=1$ and $\beta_1=3$.
Therefore, $a_1=2+\alpha_1-\beta_1=2+1-3=0$.
On the other hand, by symmetry, when $\Gamma_2=(e_{22},e_{11},e_{12},b_6,b_5,b_1,b_7,e_{24},b_4,-b_6)$, then $a_2=0$.
}
\end{rei}

\begin{center}
\includegraphics[width=9cm]{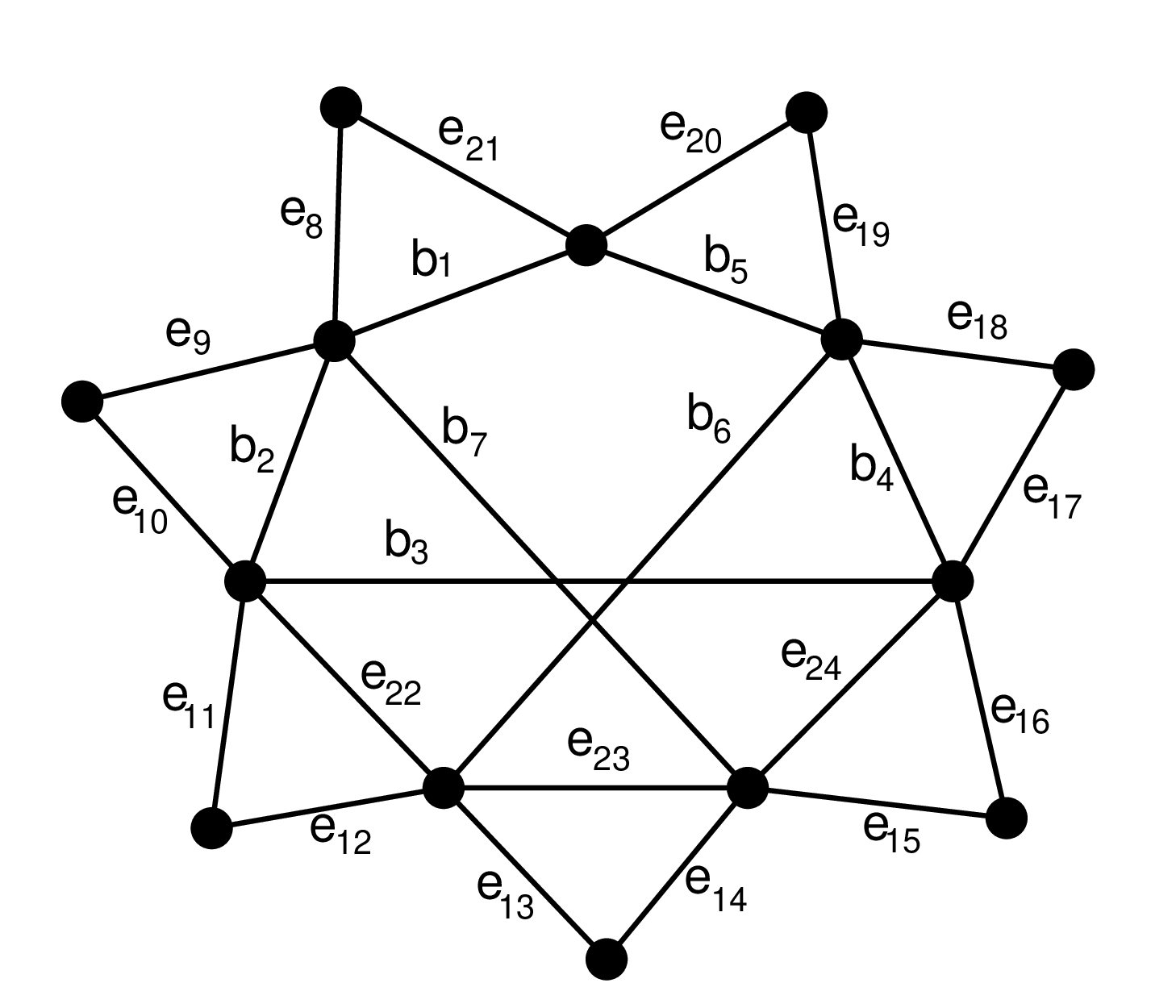}
\end{center}

\newpage
\section{The algorithm and program}
We have implemented a program for the computer algebra system Magma \cite{Magma} that determines whether a given fundamental FHM graph satisfies our criteria.
By using the program``cycle 12.c" (see ``http://sloppyjoe9.wixsite.com/mysite/program"), we tested 10 fundamental FHM graphs in appendix A.

\begin{table}[H]
\caption{Whether 10 fundamental FHM graph satisfy our criteria}
\begin{center}
\begin{tabular}{|l|r|r|r|r|}\hline
Graph number & Theorem 2.1 & Corollary 2.2 & Corollary 2.3 & Corollary 2.4 \\ \hline
Graph 1 & $\bigcirc$ & $\bigcirc$ & $\bigcirc$ & $\bigcirc$ \\ \hline
Graph 2 & $\times$ & $\times$ & $\times$ & $\times$ \\ \hline
Graph 3 & $\bigcirc$ & $\bigcirc$ & $\bigcirc$ & $\bigcirc$ \\ \hline
Graph 4 & $\bigcirc$ & $\bigcirc$ & $\bigcirc$ & $\bigcirc$ \\ \hline
Graph 5 & $\bigcirc$ & $\bigcirc$ & $\bigcirc$ & $\bigcirc$ \\ \hline
Graph 6 & $\bigcirc$ & $\bigcirc$ & $\bigcirc$ & $\bigcirc$ \\ \hline
Graph 7 & $\bigcirc$ & $\bigcirc$ & $\bigcirc$ & $\bigcirc$ \\ \hline
Graph 8 & $\bigcirc$ & $\bigcirc$ & $\bigcirc$ & $\times$ \\ \hline
Graph 9 & $\bigcirc$ & $\bigcirc$ & $\bigcirc$ & $\times$ \\ \hline
Graph 10 & $\bigcirc$ & $\bigcirc$ & $\bigcirc$ & $\bigcirc$ \\ \hline

\end{tabular}
\end{center}
\end{table}

\newpage

\appendix
\section{Example of 10 fundamental FHM graphs}

\begin{center}
\includegraphics[width=7cm]{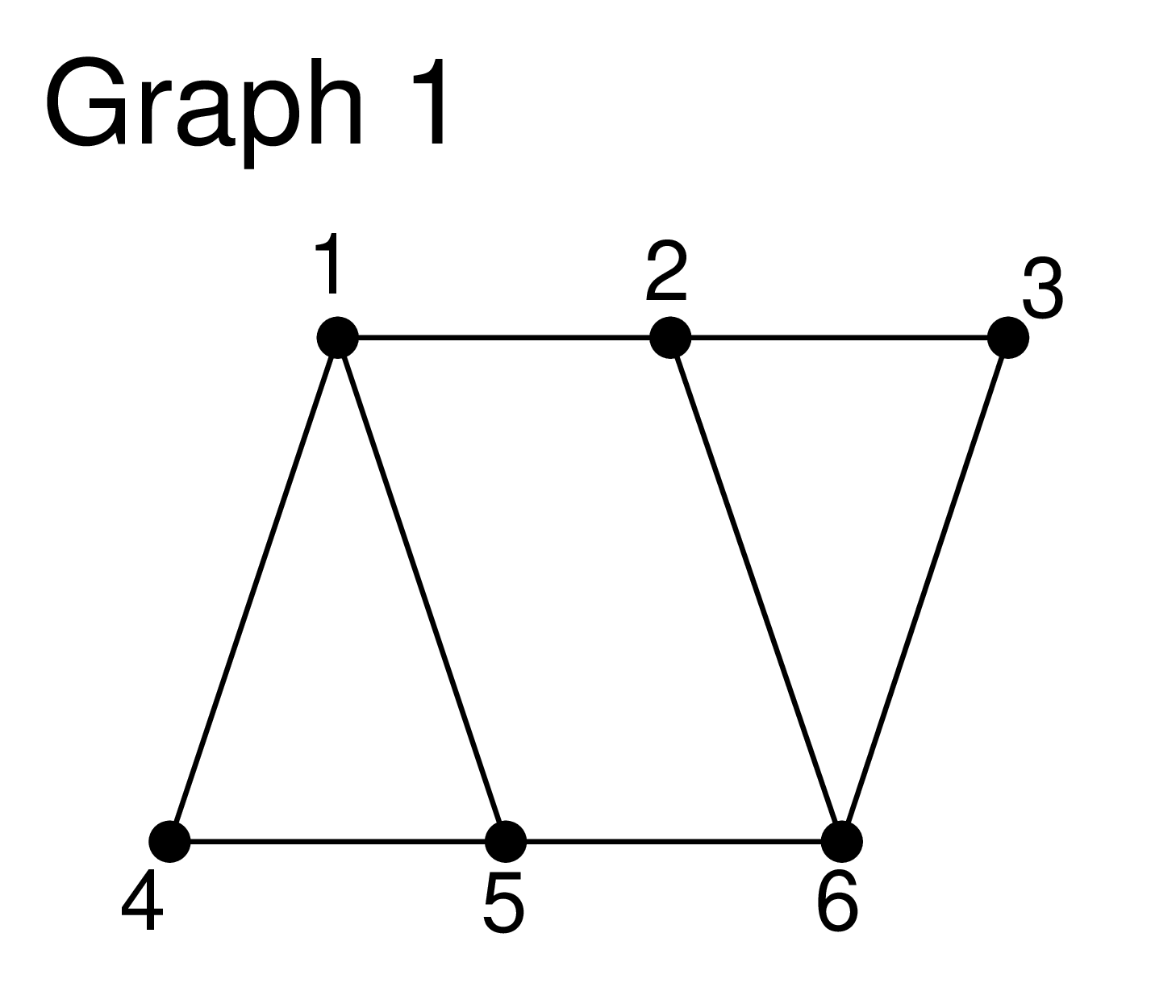} \includegraphics[width=7cm]{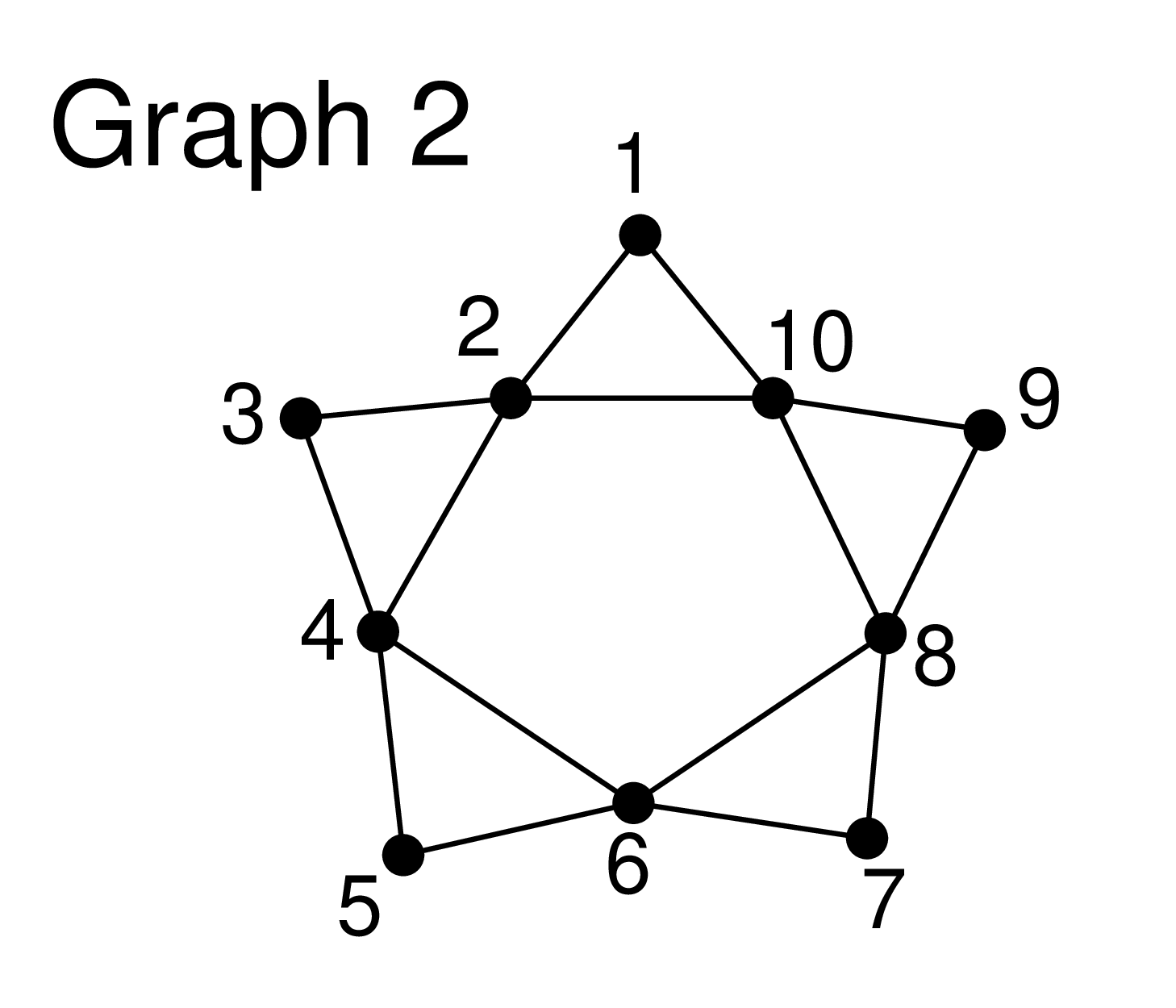}
\end{center}
\begin{center}
\includegraphics[width=7cm]{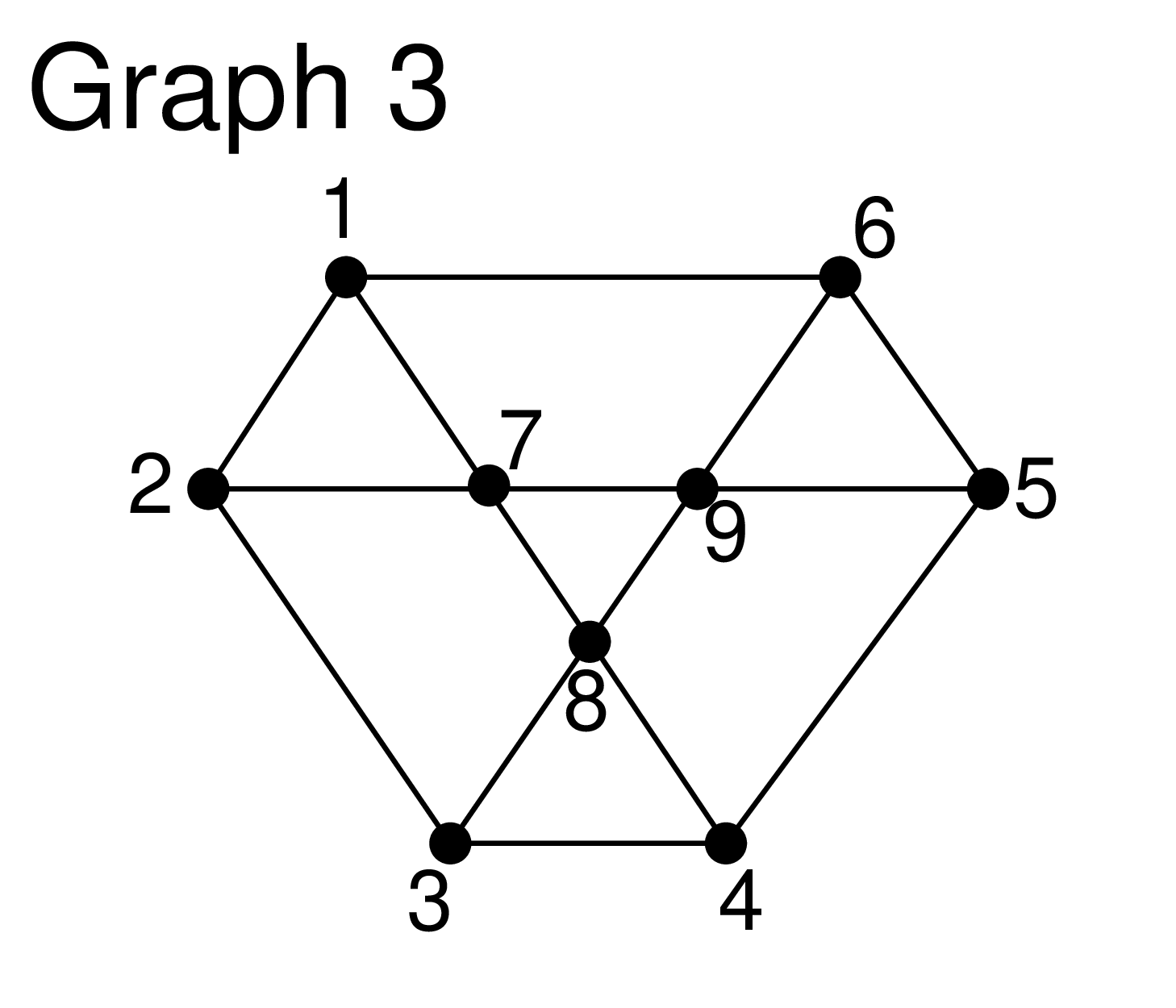} \includegraphics[width=7cm]{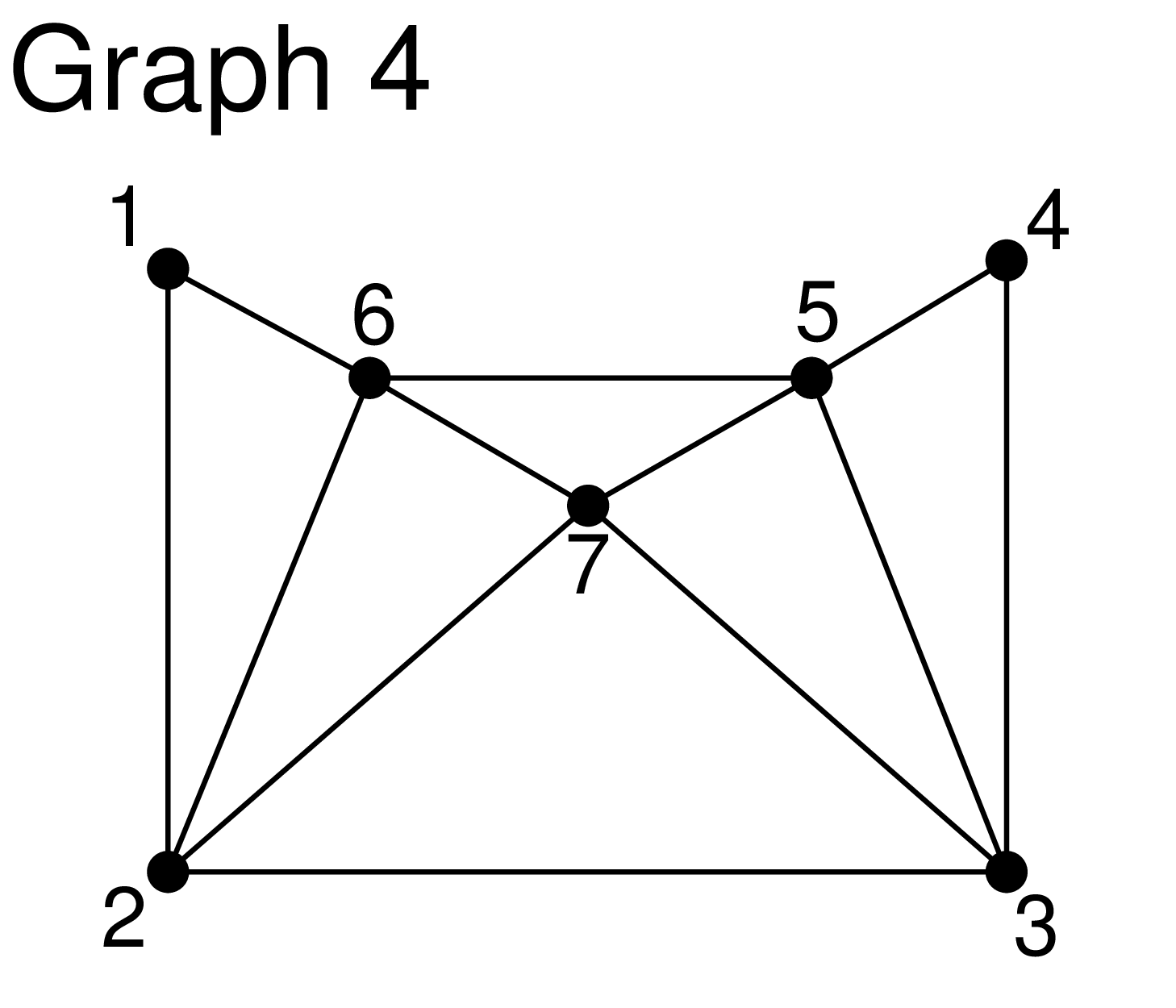}
\end{center}
\begin{center}
\includegraphics[width=7cm]{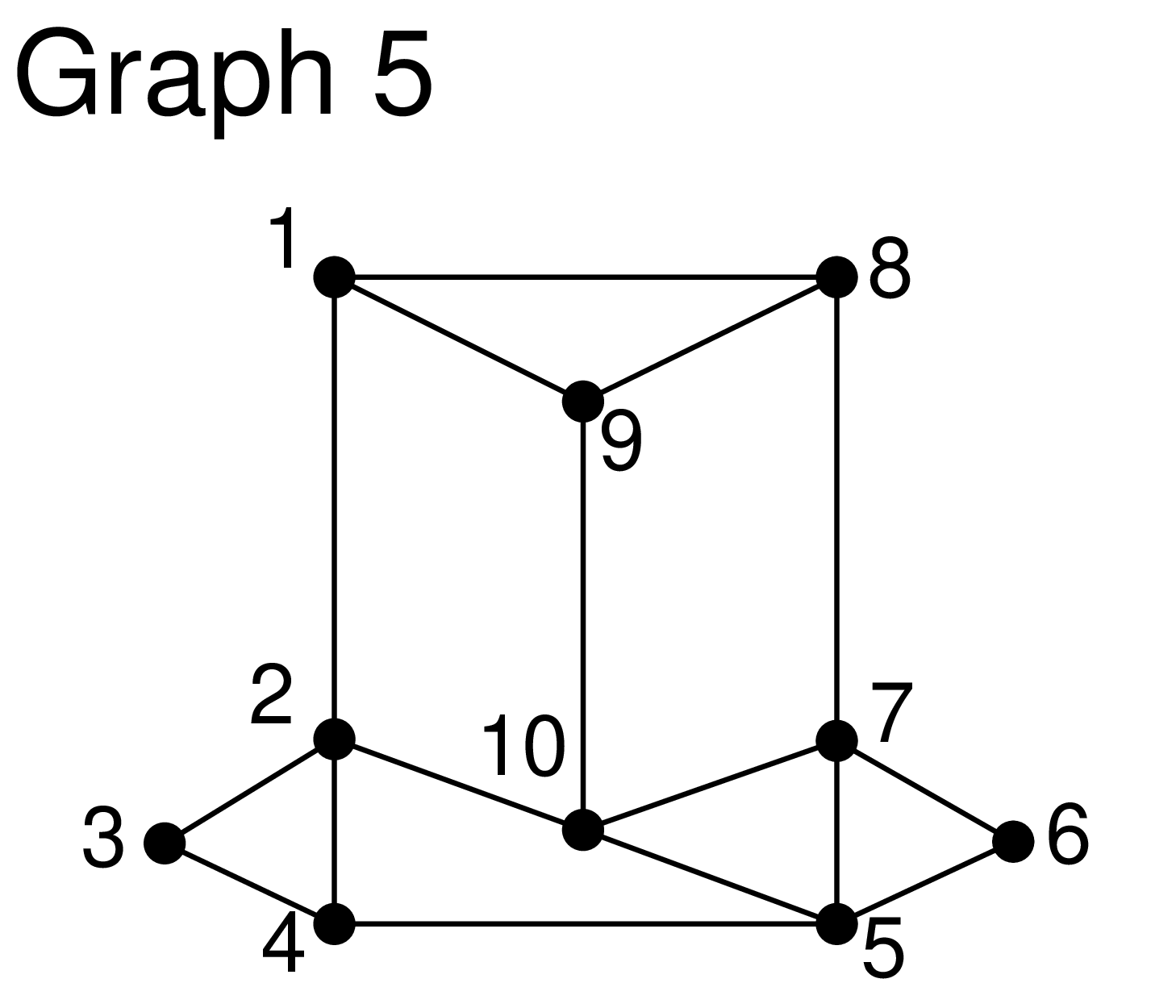} \includegraphics[width=7cm]{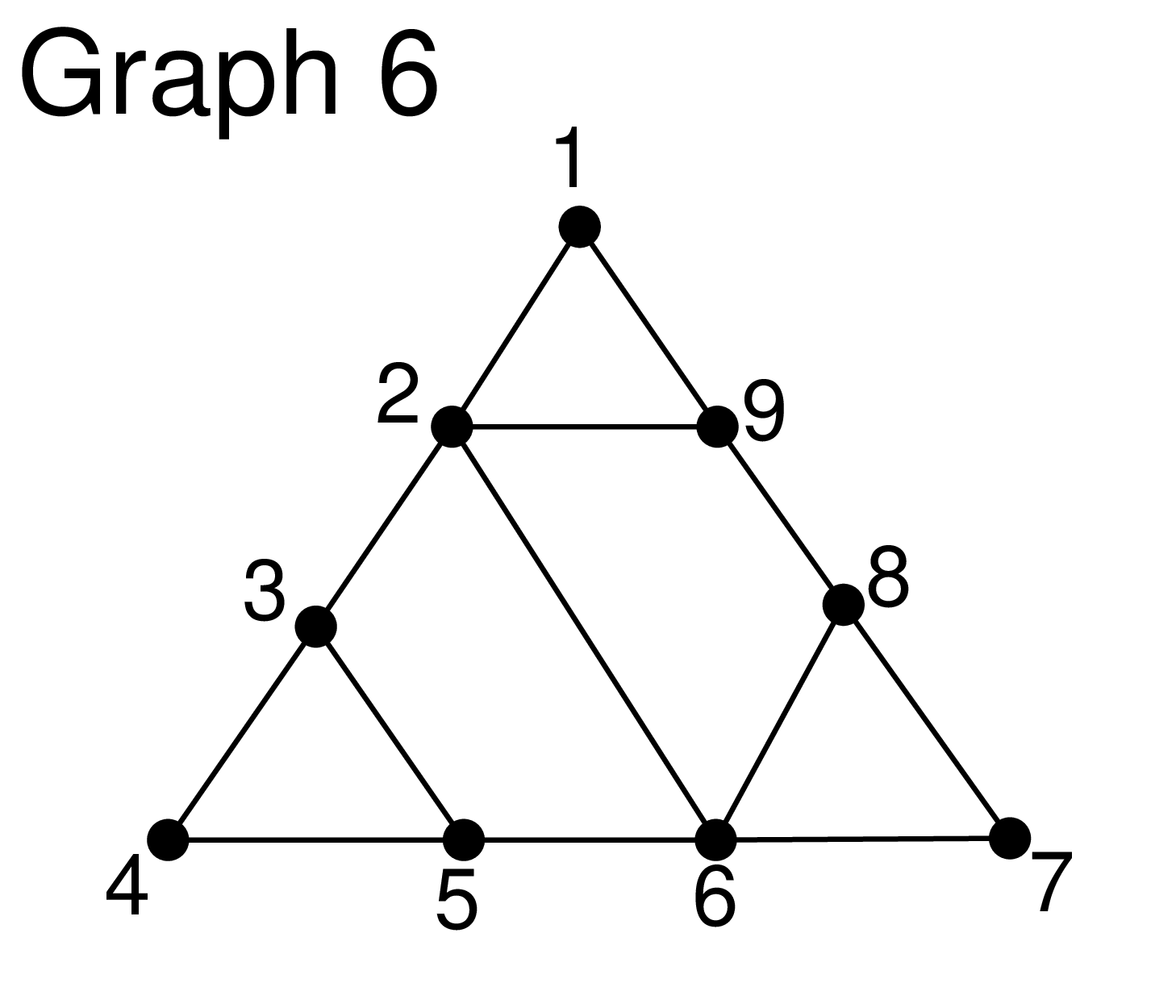}
\end{center}
\begin{center}
\includegraphics[width=7cm]{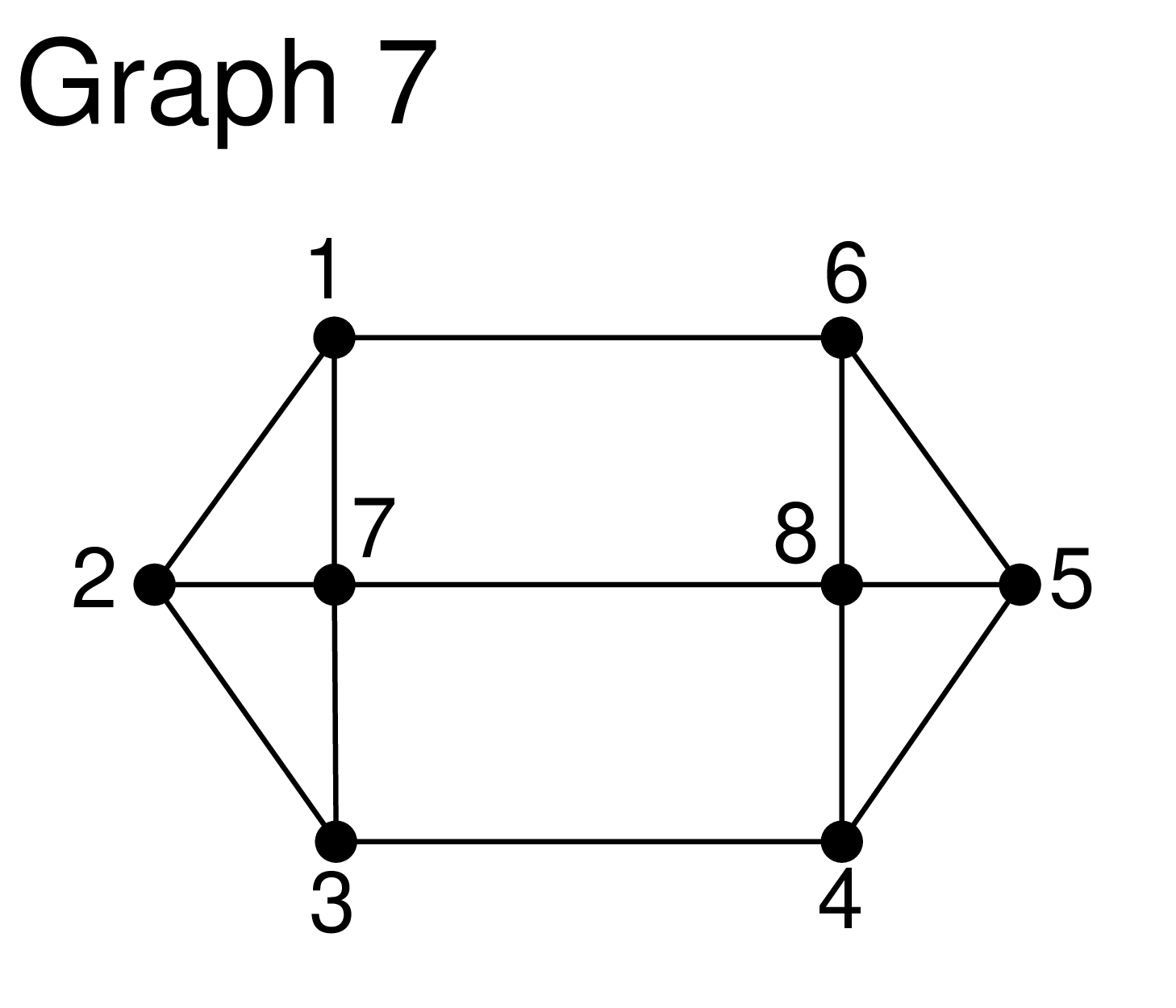} \includegraphics[width=7cm]{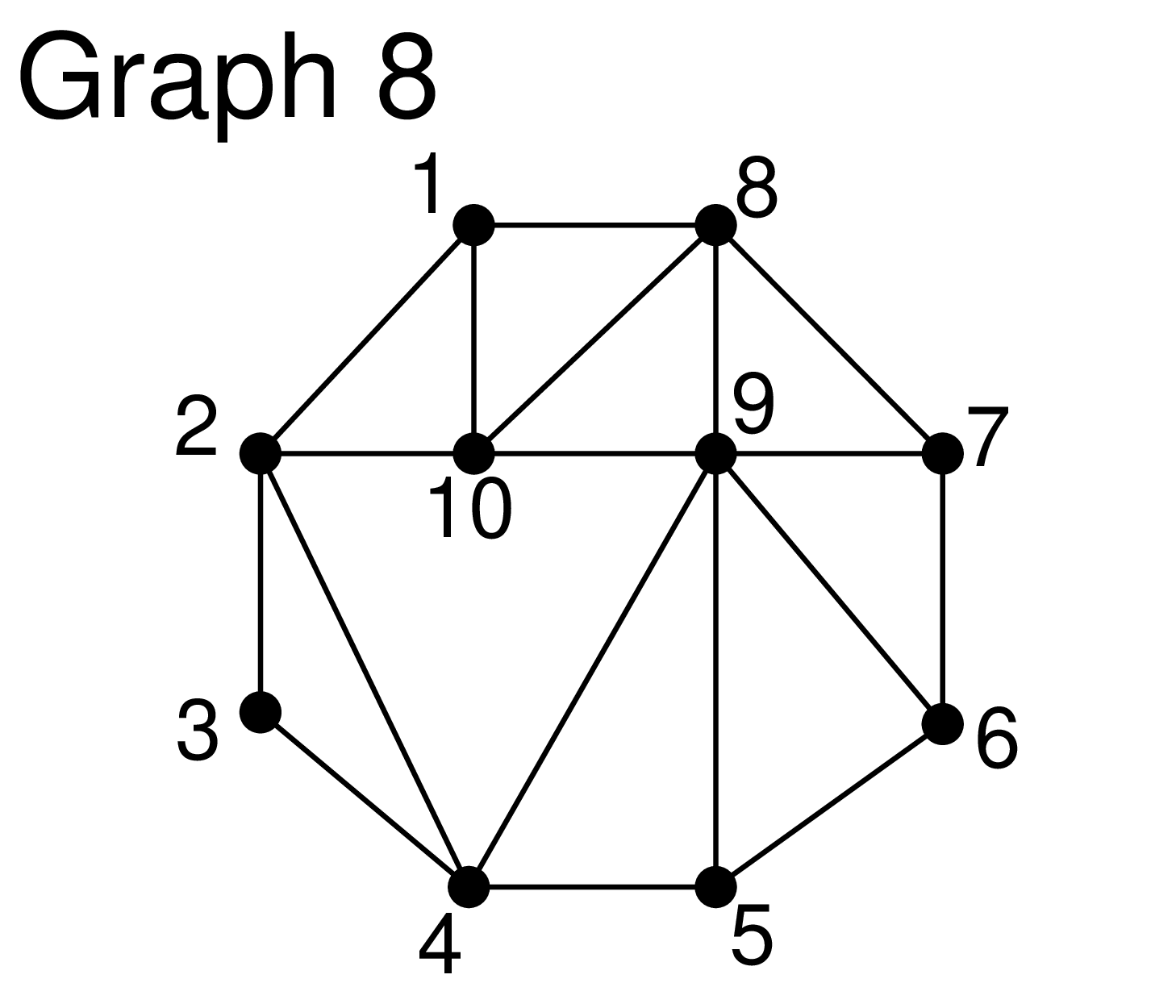}
\end{center}
\begin{center}
\includegraphics[width=7cm]{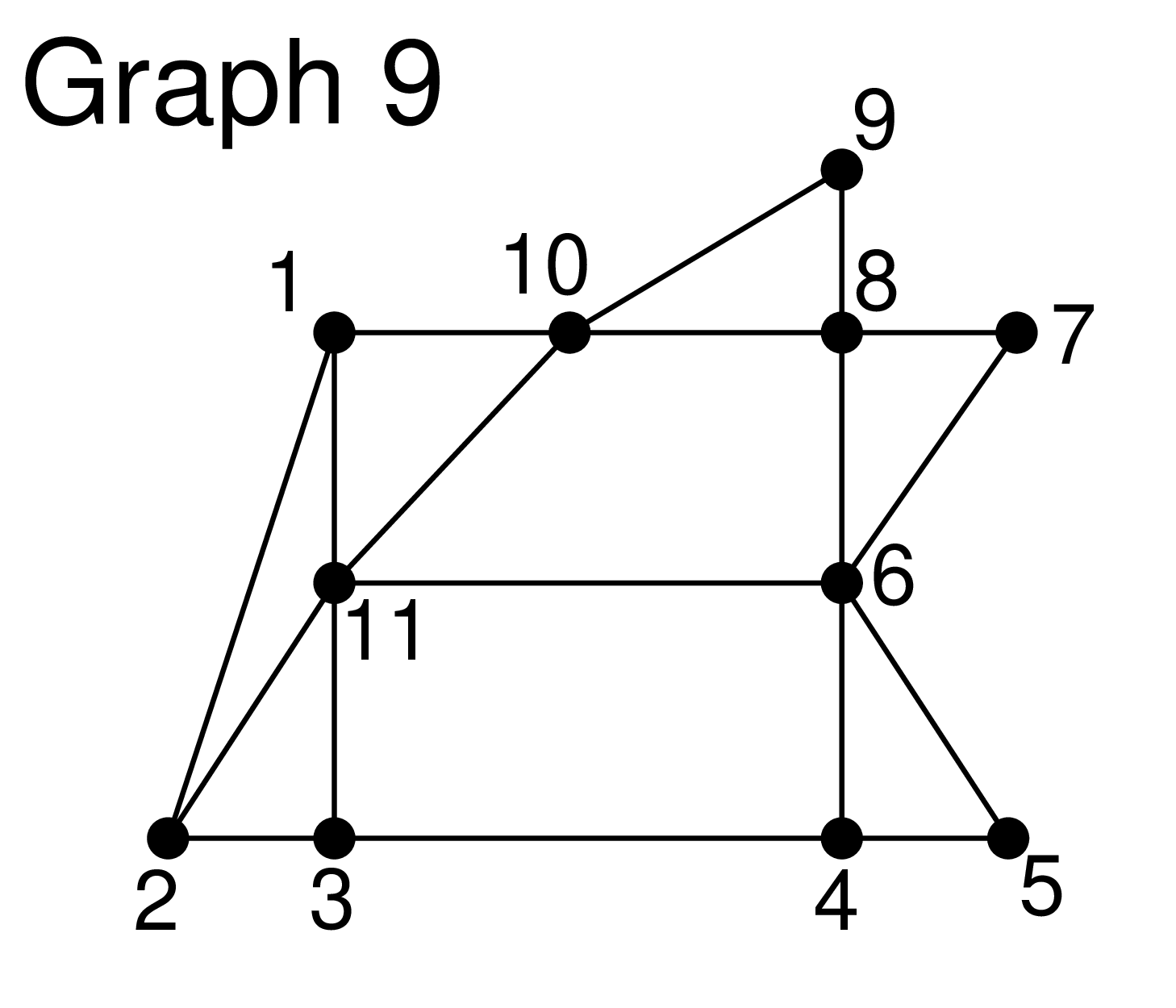} \includegraphics[width=7cm]{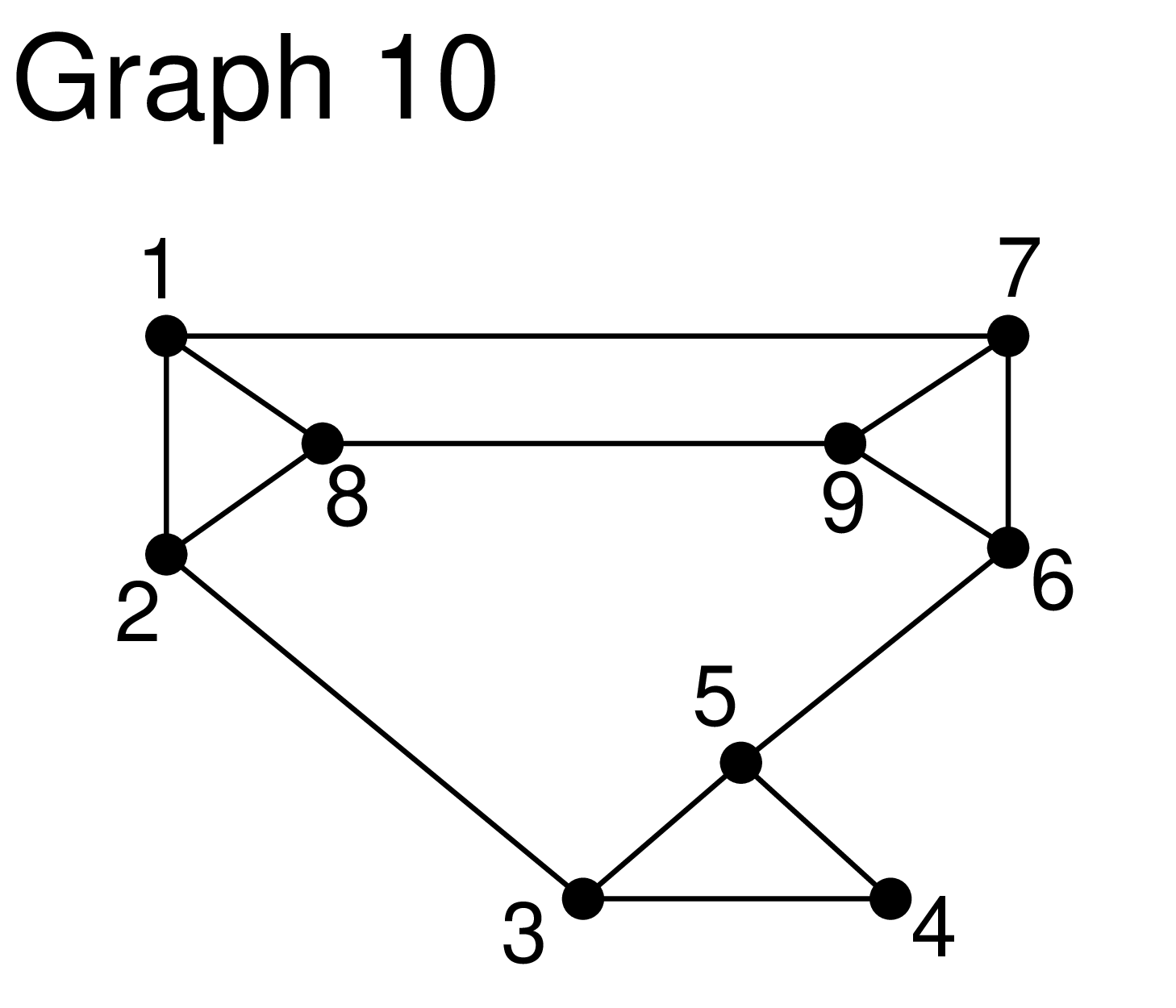}
\end{center}

(Ginji Hamano) {\sc Department of Pure and Applied Mathematics, Graduate School of Information Science and Technology, Osaka University, Suita, Osaka 565-0871, Japan }

{\it e-mail address: g-hamano@ist.osaka-u.ac.jp}\\

\end{document}